\documentclass[a4paper]{amsart}
\usepackage[ansinew]{inputenc}
\usepackage[T1]{fontenc}
\usepackage{longtable}
\usepackage{xcolor}
\usepackage{mathrsfs}
\usepackage{amssymb,amsmath,amsthm}
\usepackage[matrix, arrow, curve]{xy}
\usepackage{tikz-cd}
\usepackage{stmaryrd}
\usepackage{mathtools}
\usepackage{extarrows}
\usepackage{hhline}
\usepackage{hyperref}
\usepackage{amsmath,amssymb,tikz}
\usetikzlibrary{arrows.meta}

\definecolor{cubic}{RGB}{180,203,183}
\definecolor{quintic}{RGB}{181,199,218}
\definecolor{ghost}{gray}{0.62}
\definecolor{guide}{gray}{0.72}

\tikzset{
  quiver/.style={
    baseline=(current bounding box.center),
    x=8mm,y=8mm
  },
  vertex/.style={
    inner sep=1.2pt,
    outer sep=0pt,
    font=\normalsize
  },
  relation/.style={
    inner sep=0pt,
    outer sep=0pt,
    font=\normalsize
  },
  arr/.style={
    -{Latex[length=1.2mm,width=.8mm]},
    line width=.4pt
  },
  zero/.style={
    densely dotted,
    line width=.4pt,
    line cap=round,
    rounded corners=.8pt
  },
  longzero/.style={
    line width=3.5pt,
    line cap=round,
    line join=round
  }
}

\newcommand{\commmark}[2]{%
  \node[relation] at (#1,#2) {$\circlearrowleft$};%
}

\newcommand{\summark}[2]{%
  \begin{scope}[shift={(#1,#2)},x=1mm,y=1mm]
    \foreach \angle in {0,90,180,270}{%
      \begin{scope}[rotate=\angle]
        \draw[line width=.35pt,line cap=round]
          (.4,1.8)
          .. controls (.4,.85) and (.85,.4) ..
          (1.8,.4);
      \end{scope}
    }
  \end{scope}%
}

\usepackage{thmtools} 
\usepackage[nameinlink]{cleveref}
\usetikzlibrary{quotes,babel, angles}
\newcommand{\op}{\textnormal{op}}
\newcommand{\Hom}{\textnormal{Hom}}

\renewcommand{\mod}{\textnormal{mod}\,}

\newcommand{\im}{\textnormal{im}\,}

\newcommand{\coker}{\textnormal{coker}\,}
\newcommand{\add}{\textnormal{add}\,}

\newcommand{\cogen}{\textnormal{cogen}\,}

\numberwithin{equation}{section} \theoremstyle{plain}
\newtheorem*{thm*}{Theorem}
\newtheorem{thm}{Theorem}
\numberwithin{thm}{section}

\newtheorem*{coro*}{Corollary}

\newtheorem*{conj*}{Conjecture}
\newtheorem{lem}[thm]{Lemma}
\newtheorem*{lem*}{Lemma}
\newtheorem{prop}[thm]{Proposition}
\newtheorem*{prop*}{Proposition}
\newtheorem{rem}[thm]{Remark}
\newtheorem*{rem*}{Remark}
\newtheorem{exa}[thm]{Example}
\newtheorem*{exa*}{Example}
\newtheorem{df}[thm]{Definition}
\newtheorem*{df*}{Definition}

\newtheorem*{ques*}{Question}

\newtheorem*{construction*}{Construction}
\newtheorem*{ack*}{ACKNOWLEDGEMENTS}

\begin{document}

\title{A Solution to the Birkhoff Problem}
\author{Kevin Schlegel}
\date{}

\address{Faculty of Mathematics, Bielefeld University, 33615 Bielefeld, Germany}
\email{schlegel@math.uni-bielefeld.de}

\subjclass[2020]{16G20, 16G70, 16S90, 20K27.}
\keywords{Birkhoff Problem, monomorphism category, torsion pair}

\begin{abstract} For a commutative uniserial ring $R$ of length $n$ we describe the indecomposable objects in the monomorphism category of $R$ independent of the specific ring for $n\leq 5$ and only dependent on the residue field for $n = 6$. Thus, the classification problem for the monomorphism category of $\mathbb{Z}/p^6 \mathbb{Z}$ reduces to the known case of $\mathbb{F}_p[x]/(x^6)$, resolving the Birkhoff Problem.
\end{abstract}
\maketitle \section*{Introduction}
The study of subgroups of abelian groups goes back to the work of Miller \cite{miller, miller2} and Hilton \cite{hilton} at the beginning of the 20th century and was later developed from a more structural point of view by Birkhoff \cite{birkhoff}. In modern terms, these problems are naturally studied via monomorphism categories, see \cite{kvamme} for a detailed survey by Kvamme. The monomorphism category $\mathcal{S}(R)$ of a commutative ring $R$ has monomorphisms $X\hookrightarrow Y$ between finite length $R$-modules as objects, and suitable commutative squares as morphisms. The monomorphism category of $\mathbb{Z}$ corresponds to the setting of subgroups of finite abelian groups. Birkhoff  \cite{birkhoff} already observed that its study reduces to the monomorphism categories $\mathcal{S}(\mathbb{Z}/p^n \mathbb{Z})$ for all prime numbers $p$ and natural numbers $n$. A central question is the classification of all indecomposable objects inside $\mathcal{S}(\mathbb{Z}/p^n \mathbb{Z})$. 

Consider the general setting of a commutative uniserial ring $R$ of length $n\in \mathbb{N}$. For $n\geq 7$ Ringel and Schmidmeier showed that the monomorphism category of $R$ has wild behavior \cite{RS}, so a classification of indecomposable objects is impossible. For $n\leq 5$ Plahotnik \cite{plahotnik} proved that there are, up to isomorphism, only finitely many indecomposable objects in $\mathcal{S}(R)$. A precise classification in this finite setting for $R= \mathbb{Z}/p^n \mathbb{Z}$ and $n \leq 5$ was provided by Richman and Walker \cite{RW, RW2} based on the theory of $p$-valuated groups and valuated trees. Thus, the only unresolved case in the context of abelian groups is the monomorphism category of $\mathbb{Z}/p^6 \mathbb{Z}$. This is called the Birkhoff Problem. We solve it completely by reducing it to the known case of $R = \mathbb{F}_p[x]/(x^6)$. Ringel and Schmidmeier \cite{RS2} classified the indecomposable objects in $\mathcal{S}(k[x]/(x^n))$ for a field $k$ and $n \leq 6$. The main advantage here is that $k[x]/(x^n)$ is a finite dimensional algebra and, in particular, covering techniques are available to tackle the problem. Indeed, it was shown that the monomorphism category of $k[x]/(x^6)$ has tame behavior and the structure of the indecomposable objects can be organized via tubular families 
\cite{RS2}. Our main result is as follows.
\\
\\
\textbf{Theorem A.} (\Cref{T: main indep}) \emph{Let $R$ be a commutative uniserial ring of length $n \leq 6$. There exists a finite dimensional $k$-algebra $B$ only dependent on the residue field $k = R/(\pi)$ of $R$ that admits an equivalence between $\mathcal{S}(R)/\langle \mathcal{S}(R') \rangle$ and the category of finite dimensional torsionless $B$-modules, where $R' = R/(\pi^{n-1})$. In particular, the isomorphism classes of indecomposable objects in $\mathcal{S}(R)$ naturally correspond to those in $\mathcal{S}(k[x]/(x^n))$.} \\

The expression $\mathcal{S}(R)/\langle \mathcal{S}(R') \rangle$ in the theorem is an ideal quotient, where $\mathcal{S}(R')$ is naturally identified with a subcategory of $\mathcal{S}(R)$. A major advantage is that this category is $k$-linear and we completely describe it by a finite dimensional algebra $B$ that only depends on the residue field of $R$. This gives equivalences
\begin{align*}
    \mathcal{S}(R)/ \langle \mathcal{S}(R') \rangle \simeq \{X\in \mod B \text{ torsionless}\} \simeq \mathcal{S}(k[x]/(x^n))/ \langle \mathcal{S}(k[x]/(x^{n-1}))\rangle
\end{align*}
and so the indecomposable objects in $\mathcal{S}(R) \setminus \mathcal{S}(R')$ naturally correspond to those in $\mathcal{S}(k[x]/(x^n))\setminus \mathcal{S}(k[x]/(x^{n-1}))$. Thus, the last statement of the theorem follows inductively. Our approach is inspired by the work of Gao, K\"ulshammer, Kvamme, and Psaroudakis \cite{GKKP}, who also use ideal quotients in a similar context. In fact, they already cover the small cases $n\leq 3$.

For $n \leq 5$ we obtain a description of the finitely many indecomposable objects in $\mathcal{S}(R)$ that only uses $0$ and powers of $\pi$ as coefficients. Hence, this is completely independent of the ring $R$. For $n = 6$ the equivalences provide a description of the infinitely many indecomposable objects in $\mathcal{S}(R)$ in terms of those in $\mathcal{S}(k[x]/(x^6))$. Hence, this only depends on the residue field $k$. Moreover, it also follows that the Auslander-Reiten quiver of $\mathcal{S}(R)$ coincides with that of $\mathcal{S}(k[x]/(x^6))$ expect possibly for finitely many arrows,  see for example \cite[Lemma 1.1]{jorgensen}.

The fact that the ideal quotient in Theorem A can be described \mbox{through a finite} dimensional algebra is based on a result about torsion-free classes and certain ideal quotients established in this work. Note that the monomorphism category $\mathcal{S}(R)$ can be realized as a functorially finite torsion-free class $\mathcal{F} \subseteq \mod A$ consisting of all torsionless modules in $\mod A$, see \cite[Theorem 2.6]{kvamme}, for an Artin algebra $A$. \mbox{Moreover, under} this identification $\mathcal{S}(R') \subseteq \mathcal{S}(R)$ coincides with a full additive subcategory $\mathcal{C} \subseteq \mathcal{F}$ closed under submodules. Thus, the ideal quotient appearing in Theorem A is precisely of the form as in the following result.
\\
\\
\textbf{Theorem B.} (\Cref{T: main}) \emph{Let $A$ be an Artin algebra, $\mathcal{F} \subseteq \mod A$ a functorially finite torsion-free class and $\mathcal{C} \subseteq \mathcal{F}$ a full additive subcategory closed under submodules that contains only finitely many isomorphism classes of indecomposable modules. The ideal quotient $\mathcal{F}/ \langle \mathcal{C} \rangle$ is equivalent to a functorially finite torsion-free class $\mathcal{G} \subseteq \mod B$ for an Artin algebra $B$. Moreover, if $\mathcal{F} \subseteq \mod A$ consists of all torsionless modules, then so does $\mathcal{G} \subseteq \mod B$.}
\\

The first statement in Theorem B is based on \cite[Theorem 5.1]{enomoto} by Enomoto. He already covers the case $\mathcal{F} = \mod A$ and we extended it to an arbitrary $\mathcal{F}$ such that it can be applied to monomorphism categories. The finite dimensional algebra $B^\op$ can be explicitly calculated as the endomorphism ring of certain objects in the ideal quotient. We used AI to carry out and verify these computations, see the declaration below.
\\
\\
\textbf{Use of AI.} GPT-6 Astra (xhigh) was used to perform the calculations as described in \emph{Computations} of Section 3 and formally verify them in Lean 4 using mathlib. The source is available at \href{https://github.com/kevin-schlegel/Birkhoff-Problem}{https://github.com/kevin-schlegel/Birkhoff-Problem}. The AI was also used to visualize the quivers and relations in \Cref{R: final rem}.\\
\\
\textbf{Acknowledgments.} I am grateful to Julian K\"ulshammer for introducing me to the Birkhoff Problem and for many stimulating discussions. This work was funded by the Deutsche Forschungsgemeinschaft (DFG, German Research Foundation) under the Walter Benjamin Programme -- Project number 580964112.

\section{Torsion-free classes and ideal quotients}
We will later study monomorphism categories via their realization as torsion-free classes. In this section, we show the desired results about torsion-free classes and certain ideal quotients. 

Throughout, let $A$ be an Artin algebra (the center of $A$ is artinian and $A$ is finitely generated as a module over its center) and $\mod A$ the category of finitely generated (left) $A$-modules. A full subcategory $\mathcal{F} \subseteq \mod A$ is a \emph{torsion-free class} if $\mathcal{F}$ is closed under extensions and submodules. We are interested in the following approximation conditions.

\begin{df}\rm Let $\mathcal{C} \subseteq \mod A$ be a full additive subcategory (closed under direct sums and summands).
\begin{itemize}
    \item[(1)] A morphism $f\colon X \to C_X$ inside $\mathrm{mod}\,A$ is a \emph{left $\mathcal{C}$-approximation} of $X$ if $C_X \in \mathcal{C}$ and $\mathrm{Hom}_A(f,C)$ is surjective for all $C\in \mathcal{C}$. If every $X \in \mod A$ admits a left $\mathcal{C}$-approximation, then $\mathcal{C}$ is \emph{covariantly finite}. 
    \item[(2)] A morphism $f\colon C_X \to X$ inside $\mathrm{mod}\,A$ is a \emph{right $\mathcal{C}$-approximation} of $X$ if $C_X \in \mathcal{C}$ and $\Hom_A(C,f)$ is surjective for all $C\in \mathcal{C}$. If every $X\in \mod A$ admits a right $\mathcal{C}$-approximation, then $\mathcal{C}$ is \emph{contravariantly finite}. 
    \item[(3)] If $\mathcal{C}$ is both co- and contravariantly finite, then $\mathcal{C}$ is \emph{functorially finite}.
\end{itemize} 
\end{df}
  
The next result by Auslander and Smal{\o} \cite{ASm} shows that torsion-free classes enjoy nice approximation properties, since they are closed under submodules. 

\begin{prop} \label{P: app}\cite[Proposition 4.7 and Proposition 4.8]{ASm} Let $\mathcal{C}\subseteq \mod A$ be a full additive submodule-closed subcategory.
\begin{itemize}
    \item[(1)] Every $X\in \mod A$ admits an epic left $\mathcal{C}$-approximation $X \twoheadrightarrow  C_X$. In particular $\mathcal{C}$ is covariantly finite.
    \item[(2)] There exists $M \in \mod A$ that cogenerates $\mathcal{C}$ (so $\mathcal{C} = \cogen M$) if and only if $\mathcal{C}$ is functorially finite.
\end{itemize}
\end{prop}

\begin{exa}\label{E: ideal}\rm Let $\mathcal{C} \subseteq \mod A$ be a full additive submodule-closed subcategory. The canonical map $A \twoheadrightarrow A/I$ is a left $\mathcal{C}$-approximation where $I = \mathrm{Ann}_A(\mathcal{C})$ is the  ideal of all $a\in A$ with $aX = 0$ for all $X\in \mathcal{C}$.
\end{exa}

Auslander and Smal{\o} \cite{ASm2} also showed that a functorially finite extension closed subcategory $\mathcal{C} \subseteq \mod A$ has almost split sequences. In particular, there exists a relative Auslander-Reiten translation for $\mathcal{C}$, which we denote by $\tau_\mathcal{C}$ and its inverse by $\tau^-_\mathcal{C}$ (we write $\tau$ and $\tau^-$ for $\mathcal{C} = \mod A$). We are mainly interested in the inverse Auslander-Reiten translation $\tau^-_\mathcal{F}$ for a functorially finite torsion-free class $\mathcal{F}\subseteq \mod A$. Note that $\tau^-_\mathcal{F} X = C_{\tau^-X}$ for $X\in\mathcal{F}$ where $\tau^- X \twoheadrightarrow C_{\tau^-X}$ is the epic left $\mathcal{F}$-approximation of $\tau^- X$. The following is the main result for studying monomorphism categories later. 

\begin{thm}\label{T: main} Let $\mathcal{F} \subseteq \mod A$ be a functorially finite torsion-free class and let $\mathcal{C} \subseteq \mathcal{F}$ be a full additive submodule-closed subcategory that contains only finitely many isomorphism classes of indecomposable modules. 
\begin{itemize}
    \item[(1)] The ideal quotient $\mathcal{F}/ \langle \mathcal{C} \rangle$ is equivalent to a functorially finite torsion-free class $\mathcal{G} \subseteq \mod B$ for an Artin algebra $B$.
    \item[(2)] We may choose $B = \mathrm{End}_{\mathcal{F}/\langle \mathcal{C}\rangle}(X_1 \oplus \dots \oplus X_n)^\op$ where $X_1, \dots , X_n$ are, up to isomorphism, all indecomposable modules in $(\add A/I \cup \tau_\mathcal{F}^- \mathcal{C}) \setminus \mathcal{C}$ for the ideal $I = \mathrm{Ann}_A(\mathcal{F})$.  
    \item[(3)] With the above choice, if $\mathcal{F} = \cogen A$ then $\mathcal{G} = \cogen B$.
\end{itemize}
\end{thm}

\Cref{T: main} (1)--(2) was shown by Enomoto \cite[Theorem 5.1]{enomoto} for $\mathcal{F} = \mod A$. In the theorem $\add A/I$ denotes the collection of all direct summands of direct sums of $A/I$ and ideal quotients are defined as follows. 

For full additive subcategories $\mathcal{C}\subseteq \mathcal{D} \subseteq \mod A$ the objects in $\mathcal{D}/ \langle \mathcal{C} \rangle$ are the same as in $\mathcal{D}$ and $\mathrm{Hom}_{\mathcal{D}/\langle \mathcal{C} \rangle} (X,Y) = \mathrm{Hom}_A(X,Y)/\langle \mathcal{C} \rangle (X,Y)$ for $X,Y \in \mathcal{D}$ where $\langle \mathcal{C} \rangle (X,Y)$ denotes all morphisms $X\to Y$ factoring through a module in $\mathcal{C}$. Before we give a full proof of \Cref{T: main}, we need two lemmata about ideal quotients. The first one is folklore, see for example \cite[Lemma 1.1]{jorgensen}.

\begin{lem} \label{L: folklore} Let $\mathcal{C}\subseteq \mathcal{D}\subseteq \mod A$ be full additive subcategories.
\begin{itemize}
    \item[(1)] The quotient category $\mathcal{D}/\langle \mathcal{C} \rangle$ is linear and Hom-finite over the center of $A$. Moreover $\mathcal{D}/\langle \mathcal{C} \rangle$ has split idempotents, so it is a Krull-Schmidt category.
    \item[(2)] If $X\cong X'$ in $\mathcal{D}/\langle \mathcal{C} \rangle$ then $X\oplus C \cong X' \oplus C'$ in $\mathcal{D}$ for some $C,C' \in \mathcal{C}$.
    \item[(3)] The indecomposable objects in $\mathcal{D}/\langle \mathcal{C} \rangle$ correspond, up to isomorphism, to the indecomposable modules in $\mathcal{D} \setminus \mathcal{C}$.
    \item[(4)] The Auslander-Reiten quiver of $\mathcal{D}/\langle \mathcal{C} \rangle$ is given by the Auslander-Reiten quiver of $\mathcal{D}$ after removing the vertices corresponding to the modules in $\mathcal{C}$ along with the arrows going in and out of these vertices.
\end{itemize}
\end{lem}

The second lemma is concerned with approximations and ideal quotients. Recall that a morphism $f$ is \emph{left minimal} if $\alpha f = f$ implies that $\alpha$ is an isomorphism. Of particular interest are left minimal $\mathcal{C}$-approximations $X\to C_X$ (left approximations that are also left minimal) for a full additive subcategory $\mathcal{C} \subseteq \mod A$. Clearly, in this case $C_X$ is uniquely determined by $X$ up to isomorphism.

\begin{lem}\label{L: app} Let $\mathcal{C} \subseteq \mathcal{D} \subseteq \mod A$ be full additive subcategories, $f\colon X\to C_X$ a left $\mathcal{D}$-approximation of $X\in \mod A$ and $\bar{f}$ the corresponding morphism in $\mod A/\langle \mathcal{C} \rangle$. If $f$ is left minimal, then $\bar{f}$ is left minimal.
\end{lem}
\begin{proof} Let $\bar{\alpha} \bar{f} = \bar{f}$ for $\alpha \colon C_X \to C_X$ in $\mod A$. Then $\alpha f -f$ factors through a module $C\in \mathcal{C}$. Thus $\alpha f - f = hg$ for some $g\colon X \to C$ and $h\colon C \to C_X$. Since $\mathcal{C}\subseteq \mathcal{D}$ it follows that $g$ factors through the left $\mathcal{D}$-approximation $f$ of $X$ and so $g = g'f$ for some $g' \colon C_X \to C$. In total 
\begin{align*}
    (\alpha - hg') f = \alpha f - hg = f.
\end{align*}
Now if $f$ is left minimal, then $\alpha - hg'$ is an isomorphism. Because $\bar{h} = 0$ it follows that $\bar{\alpha}$ is an isomorphism in $\mathcal{D}/ \langle \mathcal{C}\rangle$. Hence $\bar{f}$ is left minimal in this case.
\end{proof}

\begin{proof}[Proof of \Cref{T: main}]  (1) By \cite[Theorem 5.1]{enomoto} $\mod A / \langle \mathcal{C} \rangle\simeq \mathcal{G}'$ for a functorially finite torsion-free class $\mathcal{G}' \subseteq \mod B'$ with 
\begin{align*}
B' = \mathrm{End}_{\mod A/\langle \mathcal{C} \rangle} (X_1' \oplus \dots \oplus X_m')^\op     
\end{align*}
where $X_1', \dots , X_m'$ are, up to isomorphism, all of the indecomposable modules in $(\add A \cup \tau^- \mathcal{C})\setminus \mathcal{C}$. The equivalence identifies the modules $X_1', \dots, X_m'$ with the indecomposable projective $B'$-modules. Let $\mathcal{G}\subseteq \mod B'$ be the essential image of $\mathcal{F} \subseteq \mod A$ under $L\colon \mod A \to \mod A / \langle \mathcal{C} \rangle \simeq \mathcal{G}'$. We will show that
\begin{itemize}
    \item[(i)] $\mathcal{G}$ is closed under extensions,
    \item[(ii)] $\mathcal{G}$ is closed under submodules, and
    \item[(iii)] $\mathcal{G} = \cogen L(M)$ for $M \in \mod A$ fulfilling $\mathcal{F} = \cogen M$.
\end{itemize}
Then $\mathcal{G} \subseteq \mod B'$ is a functorially finite torsion-free class by \Cref{P: app} (2) and clearly $\mathcal{F}/\langle \mathcal{C} \rangle \simeq \mathcal{G}$ via $L$. 

(i) The short exact sequences in $\mathcal{G}'$ are just the essential images of certain short exact sequences in $\mod A$, see \cite[Theorem B.2 and Proposition B.3]{enomoto}. Thus, a short exact sequence $0 \to L(X) \to L(Y) \to L(Z) \to 0 $ in $\mathcal{G}'$ corresponds to a short exact sequence $0 \to X \oplus C_X \to Y \oplus C_Y \to Z \oplus C_Z \to 0 $ in $\mod A$ with $C_X, C_Y, C_Z \in \mathcal{C}$ by \Cref{L: folklore} (2). Now if $X,Z \in \mathcal{F}$ then also $Y \in \mathcal{F}$ since $\mathcal{F}$ is extension-closed and $\mathcal{C} \subseteq \mathcal{F}$. It follows that $\mathcal{G}$ is extension-closed. 

(ii) Let $f\colon X\to Y$ be a morphism in $\mod A$ such that $L(f)$ is monic and $Y\in \mathcal{F}$. For $\mathcal{G}$ to be submodule-closed, we must show that $X\in \mathcal{F}$. Consider the canonical inclusion $\iota \colon \ker f \hookrightarrow X$. Then $L(f)L(\iota) = 0$. Since $L(f)$ is monic, already $L(\iota) = 0$ so $\iota$ factors through some $C\in \mathcal{C}$. Because $\mathcal{F}$ is submodule-closed and $\mathcal{C}\subseteq \mathcal{F}$ it follows that $\ker f \in \mathcal{F}$. Now the short exact sequence $0 \to \ker f \to X \to \im f \to 0$ implies that $X \in \mathcal{F}$ as $\mathcal{F}$ is extension-closed and $Y \supseteq \im f \in \mathcal{F}$.

(iii) By assumption $\mathcal{F}$ is functorially finite, so $\mathcal{F} = \cogen M$ for $M\in \mod A$ by \Cref{P: app} (2). Since $\mod A \to \mod A /\langle \mathcal{C} \rangle$ preserves monomorphisms by the dual of \cite[Proposition 5.3 (3)]{enomoto}, it follows that $\mathcal{G} = \cogen L(M)$.

(2) We have established that $\mathcal{F}/\langle \mathcal{C} \rangle \simeq \mathcal{G} \subseteq \mod B'$ where $\mathcal{G}$ is a functorially finite torsion-free class. Let $J = \mathrm{Ann}_{B'}(\mathcal{G})$ and set $B = B'/J$. Clearly $\mathcal{G} \subseteq \mod B$ remains a functorially finite torsion-free class. The quotient $B' \twoheadrightarrow B$ is an epic left $\mathcal{G}$-approximation, see \Cref{E: ideal}, so it is left minimal. Note that, since $L$ is full, a left $\mathcal{F}$-approximation of $X\in \mod A$ turns into a left $\mathcal{G}$-approximation of $L(X)$ and left minimality is preserved by \Cref{L: app}. Thus, to determine $\add B$ we may consider the preimage of $\add B'$ under $L$, compute left minimal $\mathcal{F}$-approximations of modules in the preimage, and then apply $L$. Now the preimage is given by $\add (A\cup \mathcal{C} \cup \tau^- \mathcal{C})$ and $A\to A/I, \tau^- X \to \tau^-_\mathcal{F} X, 1\colon X \to X$ for $X\in \mathcal{C}$ are left minimal $\mathcal{F}$-approximations. Since $X_1, \dots , X_n$ are, up to isomorphism, all of the indecomposable modules in $(\add A/I \cup \tau_\mathcal{F}^- \mathcal{C}) \setminus \mathcal{C}$ it follows with \Cref{L: folklore} (3) that $L X_1, \dots, LX_n$ are, up to isomorphism, all of the indecomposable modules in $\add B$. Hence $B \cong \mathrm{End}_{\mathcal{F}/\langle \mathcal{C} \rangle} (X_1\oplus \dots \oplus X_n)^\op$ as $B$ is basic.

(3) Since $B \in \mathcal{G}$ we have $\cogen B \subseteq \mathcal{G}$. Suppose that $\mathcal{F} = \cogen A$. By (iii) it follows that $\mathcal{G} = \cogen L(A)$ and in (2) it was proven that $L(A/I) \in \add B$. Clearly $I = 0$ because $A \in \mathcal{F}$ and so $\mathcal{G} = \cogen B$.  
\end{proof}

\section{Application to monomorphism categories}

In this section, we study monomorphism categories of certain commutative rings using \Cref{T: main}. We recall the definition and give some motivating examples.

\begin{df} \rm The \emph{monomorphism category}, denoted by $\mathcal{S}(R)$, of a commutative ring $R$ has monomorphisms $X\hookrightarrow Y$ between finite length $R$-modules as objects. Morphisms from $X\hookrightarrow Y$ to $X' \hookrightarrow Y'$ in $\mathcal{S}(R)$ are commutative diagrams
\begin{equation*}
    \begin{tikzcd}
        X \arrow[r, hook] \arrow[d] & Y \arrow[d] \\
        X' \arrow[r, hook] & Y'
    \end{tikzcd}
\end{equation*}
of $R$-module homomorphisms.
\end{df}

\begin{exa}\label{E: main ex}\rm \begin{itemize}
    \item[(1)] The examples (2) and (3) illustrate the importance of the monomorphism category of a commutative uniserial ring $R$ of finite length. Let $n$ be the length of $R$. For $n\leq 5$ Plahotnik \cite{plahotnik} showed that $\mathcal{S}(R)$ contains, up to isomorphism, only finitely many indecomposable objects. For $n \geq 7$ Ringel and Schmidmeier  \cite{RS} proved that $\mathcal{S}(R)$ has wild behavior, so a classification of indecomposable objects is impossible.
    
    \item[(2)] The monomorphism category $\mathcal{S}(\mathbb{Z})$ gives the setting to study subgroups of finite abelian groups categorically, which was first pursued by Birkhoff \cite{birkhoff}. He observed that an indecomposable object already belongs to $\mathcal{S}(\mathbb{Z}/p^n \mathbb{Z})$ for a prime number $p$ and natural number $n$. Richman and Walker \cite{RW, RW2} classified the indecomposable objects in $\mathcal{S}(\mathbb{Z}/p^n \mathbb{Z})$ for $n \leq 5$ based on the theory of $p$-valuated groups and valuated trees. The remaining case $n = 6$ (by wildness $n\geq 7$ is not possible) is called the \emph{Birkhoff Problem}.
    \item[(3)] The monomorphism category $\mathcal{S}(k[x])$ gives the categorical framework to study square matrices with an invariant subspace. It is not hard to see that an indecomposable object in $\mathcal{S}(k[x])$ already belongs to $\mathcal{S}(k[x]/(f^n))$ for an irreducible polynomial $f\in k[x]$ and natural number $n$. We focus on the case $f = x$ studied by Ringel and Schmidmeier \cite{RS2}. (In fact, this covers various cases as $k[x]/(f^n) \cong F[y]/(y^n)$ if $ F =k[x]/(f)$ is separable over $k$.) They completely described the indecomposable objects inside $\mathcal{S}(k[x]/(x^n))$ for $n \leq 6$ based on covering techniques from representation theory of finite dimensional algebras, thereby solving the classification problem.
\end{itemize}
\end{exa}
Our goal will be to reduce the classification problem of \Cref{E: main ex} (1) to that of (3) for $n \leq 6$. In particular, this will give a solution to the Birkhoff Problem. To achieve this, we utilize \Cref{T: main}. Throughout, fix a commutative uniserial ring $R$ of finite length. 

The first step is to realize $\mathcal{S}(R)$ as a torsion-free class inside the module category of an Artin algebra. For this, consider the matrix algebra $A = \left(\begin{smallmatrix} R & R \\ 0 & R \end{smallmatrix}\right)$. Clearly $A$ is an Artin algebra, a (left) $A$-module corresponds to a morphism between $R$-modules and an $A$-module homomorphism is given by a suitable commutative diagram of $R$-module homomorphisms. 

\begin{lem}\label{L: ff} The monomorphism category of $R$ is equivalent to a functorially finite torsion-free class $\mathcal{F}\subseteq \mod A$. Moreover $\mathcal{F} = \cogen A$.
\end{lem}
\begin{proof} The first claim is a special case of \cite[Theorem 2.6]{kvamme}. Let $Q$ be a cogenerator of $\mod R$. Consider the $A$-module $M$ corresponding to $(0 \to Q) \oplus ( Q \xrightarrow{1} Q) \in \mathcal{S}(R)$. By \cite[Theorem 3.1 and Corollary 3.2]{LZ} $\mathcal{F} = \cogen M$. Since $R$ is self-injective we may choose $Q = R$. In this case $M =A$. 
\end{proof}

Next, we consider a suitable ideal quotient of the functorially finite torsion-free class $\mathcal{F} \subseteq \mod A$ associated to $\mathcal{S}(R)$. Note that $\mathcal{F}$ is an $R$-linear category, since $R$ is commutative. Let $(\pi) \subseteq R$ be a maximal ideal, $n\in \mathbb{N}$ minimal with $\pi^n =0$, and $k = R/(\pi)$ the residue field of $R$. 

\begin{lem}\label{L: sub} Let $\mathcal{C} \subseteq \mathcal{F}$ be the full subcategory of all $X\in \mathcal{F}$ with $\pi^{n-1}X = 0$.
\begin{itemize}
    \item[(1)] Then $\mathcal{C}$ is submodule-closed in $\mod A$.
    \item[(2)] The ideal quotient $\mathcal{F}/\langle \mathcal{C} \rangle$ is a $k$-linear category.
    \item[(3)] We have $\mathcal{C} \simeq \mathcal{S}(R')$ for $R' = R/(\pi^{n-1})$.
\end{itemize}
\end{lem}
\begin{proof} (1) For $X\in \mod A$ with $\pi^{n-1} X =0$ clearly $\pi^{n-1} Y =0$ for all submodules $Y\subseteq X$. Now since $\mathcal{F}$ is submodule-closed, also $\mathcal{C}$ is submodule-closed.

(2) For a morphism $f\colon X\to Y$ inside $\mathcal{F}$ we have $\pi^{n-1} \im \pi f =\im  \pi^n f = 0$. Hence $\im \pi f \in \mathcal{C}$ and $\pi f\in \langle \mathcal{C} \rangle$. It follows that  $\mathcal{F}/\langle \mathcal{C} \rangle$ is $k$-linear. 

(3) Clearly $\mathcal{S}(R')$ identifies with the full subcategory of all $X\in \mathcal{S}(R)$ fulfilling $\pi^{n-1} X =0$. Thus, the equivalence $\mathcal{S}(R) \simeq \mathcal{F}$ induces $\mathcal{S}(R') \simeq \mathcal{C}$.
\end{proof}

From now on, simply identify $\mathcal{S}(R) = \mathcal{F}\subseteq \mod A$ and $\mathcal{S}(R') = \mathcal{C} \subseteq \mathcal{F}$. We are ready to apply \Cref{T: main} to show the following result. It serves as the key ingredient to study monomorphism categories inductively. 

\begin{prop}\label{P: key} Suppose that $\mathcal{S}(R')$ contains only finitely many isomorphism classes of indecomposable objects. The ideal quotient $\mathcal{S}(R)/\langle \mathcal{S}(R')\rangle$ is equivalent to $\cogen B \subseteq \mod B$ for the finite dimensional $k$-algebra 
\begin{align*}
B = \mathrm{End}_{\mathcal{S}(R)/\langle \mathcal{S}(R') \rangle}(Y_1\oplus  \dots \oplus Y_s)^\op    
\end{align*}
where $Y_1, \dots, Y_s$ are, up to isomorphism, all indecomposable objects in
\begin{align*}
\{ 0\to R,\,R\xrightarrow{1} R\} \cup \tau_{\mathcal{S}(R)}^- \mathcal{S}(R') \setminus \mathcal{S}(R').
\end{align*}
\end{prop}
\begin{proof} Recall that we identified $\mathcal{S}(R) = \mathcal{F} \subseteq \mod A$ and $\mathcal{S}(R') = \mathcal{C} \subseteq \mathcal{F}$, and
\begin{align*}
A \cong (0\to R) \oplus (R \xrightarrow{1} R )   
\end{align*}
under this identification. By \Cref{L: ff} $\mathcal{F}$ is a functorially finite torsion-free class and $\mathcal{F} = \cogen A$. Moreover $\mathcal{C}$ is closed under submodules by \Cref{L: sub} (1) and contains only finitely many isomorphism classes of indecomposable modules by assumption. Thus, all conditions are fulfilled to apply \Cref{T: main} and it follows that $\mathcal{F}/\langle \mathcal{C} \rangle \simeq \cogen B$ for the Artin algebra $B = \mathrm{End}_{\mathcal{F}/\langle \mathcal{C} \rangle} (Y_1\oplus \dots \oplus Y_s)^\op$ where $Y_1, \dots, Y_s$ are, up to isomorphism, all of the indecomposable modules in $\add A \cup \tau_\mathcal{F}^-\mathcal{C} \setminus \mathcal{C}$. By \Cref{L: sub} (2) the ideal quotient $\mathcal{F}/ \langle \mathcal{C} \rangle$ is $k$-linear and so $B$ is a finite dimensional $k$-algebra.
\end{proof}

We will make use of the above result to completely classify the indecomposable objects inside the monomorphism category $\mathcal{S}(R)$ for $n \leq 6$ by reducing it to the known case of $\mathcal{S}(k[x]/(x^n))$, see \Cref{E: main ex}.

\begin{thm}\label{T: main indep} For $n \leq 6$ there is a finite dimensional $k$-algebra $B$ only dependent on the residue field $k$ of $R$ that admits an equivalence $\mathcal{S}(R)/\langle \mathcal{S}(R') \rangle \simeq \cogen B$. In particular, the isomorphism classes of indecomposable objects in $\mathcal{S}(R)$ naturally correspond to those in $\mathcal{S}(k[x]/(x^n))$. 
\end{thm}

The following proof of the theorem depends on computations, listed as C1--C3, that still need to be verified. This is done in Section 3.

\begin{proof} For $m \leq 6$ let $s_m = 2+2(m-1) \frac{6}{6-m}$ ($=\infty$ for $m = 6$) be the number of isomorphism classes of indecomposable objects in $\mathcal{S}(k[x]/(x^m))$, see \cite{RS2}. Suppose that we have listed, up to isomorphism, all of the indecomposable objects in $\mathcal{S}(R')$ and there are $s_{n-1}$-many of them. Denote them by $X_1, X_2, \dots, X_{s_{n-1}} \in \mathcal{S}(R')$. We will do the following computations.
\begin{itemize}
    \item[(C1)] We compute $\tau^-_{\mathcal{S}(R)} X_i$ for $i = 1, \dots, s_{n-1}$ and let $Y_1, \dots, Y_s$ be all of these indecomposable objects that are not contained in $\mathcal{S}(R')$ together with the objects $0 \to R$ and $ R \xrightarrow{1} R$.
    \item[(C2)] We give a bound path algebra $E = kQ/I$, where $Q$ and generators of $I$ are independent of the specific ring $R$, and compute an isomorphism
    \begin{align*}
            E \cong \mathrm{End}_{\mathcal{S}(R)/ \langle \mathcal{S}(R')\rangle} (Y_1 \oplus \dots \oplus Y_s). 
    \end{align*}
    \item[(C3)] If $n \neq 6$ we list $(s_n-s_{n-1})$-many non-isomorphic objects $X_{s_{n-1}+1}, \dots, X_{s_n}$ in $\mathcal{S}(R) \setminus \mathcal{S}(R')$ and compute that $\mathrm{End}_{\mathcal{S}(R)} (X_j)$ is local for $j > s_{n-1}$.
\end{itemize}

Suppose that C1--C3 have been achieved. Then by \Cref{P: key} we obtain an equivalence 
\begin{align*}
    \mathcal{S}(R)/ \langle \mathcal{S}(R') \rangle \simeq \cogen E^\op \simeq \mathcal{S}(k[x]/(x^n))/ \langle \mathcal{S}(k[x]/(x^{n-1})) \rangle
\end{align*}
as $Q$ and generators of $I$ are independent of $R$. By \Cref{L: folklore} (3) the number of isomorphism classes of indecomposable objects in $\mathcal{S}(R) \setminus \mathcal{S}(R')$ equals $s_n-s_{n-1}$. Thus $X_1, \dots, X_{s_n}$ are, up to isomorphism, all indecomposable objects inside $\mathcal{S}(R)$. The theorem therefore follows inductively over $n$.  
\end{proof}

\begin{rem}\rm By \Cref{T: main indep} the indecomposable objects in $\mathcal{S}(\mathbb{Z}/p^6 \mathbb{Z})$ naturally correspond to those of $\mathcal{S}(\mathbb{F}_p[x]/(x^6))$ providing a solution to the Birkhoff Problem. Moreover, the equivalence 
\begin{align*}
    \mathcal{S}(\mathbb{Z}/p^6 \mathbb{Z})/ \langle \mathcal{S}(\mathbb{Z}/p^5 \mathbb{Z}) \rangle \simeq \cogen E^\op \simeq \mathcal{S}(\mathbb{F}_p[x]/(x^6))/ \langle \mathcal{S}(\mathbb{F}_p[x]/(x^{5})) \rangle
\end{align*}
shows that the Auslander-Reiten quiver of $\mathcal{S}(\mathbb{Z}/p^6 \mathbb{Z})$ coincides with the Auslander-Reiten quiver of $\mathcal{S}(\mathbb{F}_p[x]/(x^6))$ up to finitely many arrows, see \Cref{L: folklore} (4). The latter was determined by Ringel and Schmidmeier \cite{RS2}.
\end{rem}

\section{Computations}

Our goal is to carry out the computations C1--C3 in the proof of \Cref{T: main indep}. As before, let $R$ be a commutative uniserial ring of finite length. Moreover, let $(\pi) \subseteq R$ be a maximal ideal, $n\in \mathbb{N}$ minimal with $\pi^n = 0$, and $k = R/(\pi)$ the residue field of $R$.  

We use the following notation to describe the objects inside the monomorphism category of $R$. Set $C_i = R/(\pi^{i})$ for $i \leq n$. A morphism $f\colon C_i\to C_j$ is uniquely determined by $f(\bar{1}) = \bar{r}$ for $r\in R$ and we simply identify $f = r$. In that way, an arbitrary morphism 
\begin{align*}
f\colon \bigoplus_{i=1}^r  C_{a_i} \longrightarrow \bigoplus_{j=1}^s C_{b_j}    
\end{align*}
is given by an $s \times r$ matrix over $R$. Note that, up to isomorphism, every object in $\mathcal{S}(R)$ can be described in this way. 

\begin{exa}\label{E: C3} \rm The Auslander-Reiten quiver of $\mathcal{S}(k[x]/(x^5))$ was determined by Ringel and Schmidmeier \cite{RS2}. It is given by the following depiction from \cite{RS2}.
\begin{center}
\includegraphics[scale=0.4]{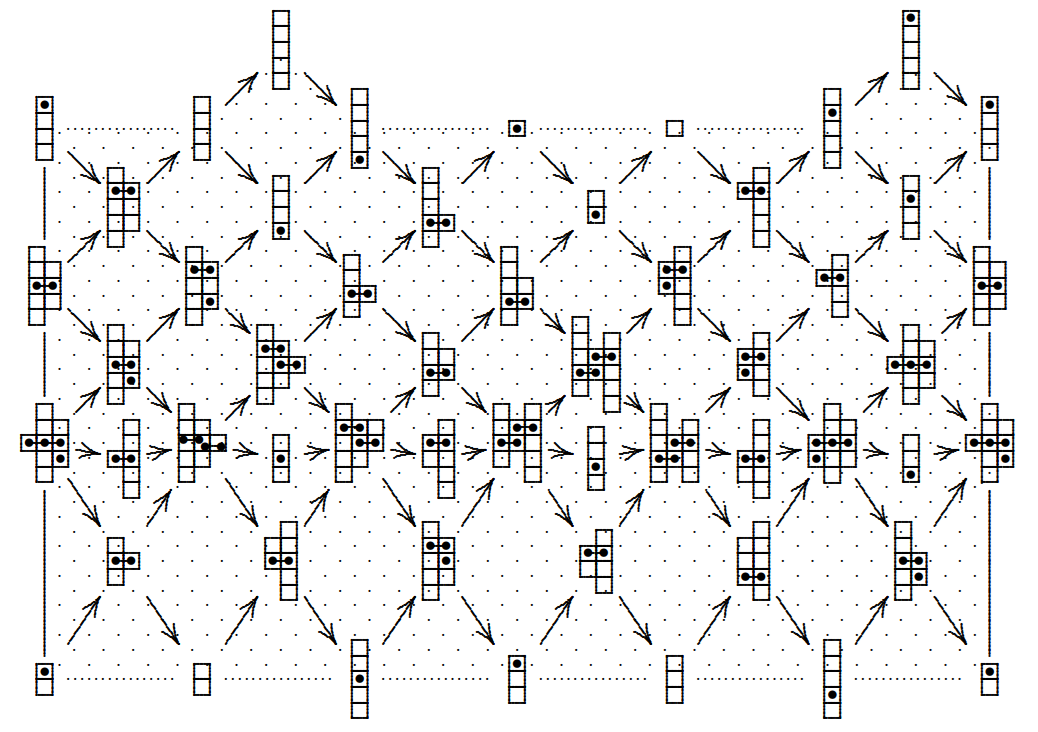}
\end{center}

In particular, we can read off all of the indecomposable objects in $\mathcal{S}(k[x]/(x^5))$. The necessary translation is clear from the examples
\begin{equation*}
    \vcenter{\hbox{\includegraphics[scale=0.4]{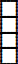}}} \,\, \hat{=} \,\,  0 \xrightarrow{0} C_4, \quad \vcenter{\hbox{\includegraphics[scale=0.4]{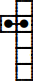}}}\,\,\hat{=} \,\,  C_4 \xrightarrow{\left(\begin{smallmatrix}
        1\\
        x
    \end{smallmatrix}\right)} C_1 \oplus C_5, \quad    
    \vcenter{\hbox{\includegraphics[scale=0.4]{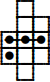}}} \,\, \hat{=} \,\, C_1 \oplus C_3 \xrightarrow{\left(\begin{smallmatrix}
        x & 1\\
        0 & x^2\\
        0 & x
    \end{smallmatrix}\right)} C_2 \oplus C_5 \oplus C_3.
\end{equation*}
 The notation still makes sense to describe objects in $\mathcal{S}(R)$ after replacing $x^i$ with $\pi^i$ (as long as the height of the figure does not exceed $n$). This will be used to list the indecomposable objects for the computation C3.   
\end{exa}

The following lemma will be used for the computation C2.

\begin{lem} \label{L: C2} A morphism
\begin{equation*}
    \begin{tikzcd}
        X \arrow[r, hook] \arrow[d, "\alpha", swap] & Y \arrow[d, "\beta"] \\
        X' \arrow[r, hook] & Y' 
    \end{tikzcd}
\end{equation*}
in $\mathcal{S}(R)$ is contained in $\langle \mathcal{S}(R') \rangle$ if and only if $\pi^{n-1} \beta = 0$.
\end{lem}
\begin{proof} The only if part is clear as every $R'$-module is annihilated by $\pi^{n-1}$. For the other implication, observe that the square factors through $f\colon \im \alpha \to \im \beta $. Now $f$ is always a monomorphism that is contained in $\mathcal{S}(R')$ if $\pi^{n-1} \beta = 0$.
\end{proof}

Lastly, we want to simplify the computation C1.

\begin{lem}\label{L: tau} Let $f\colon X\to Y$ be an indecomposable object inside $\mathcal{S}(R')$ and write $\coker f = \bigoplus_{i=1}^r C_{a_i}$. For each $i$ choose $y_i\in Y$ lifting $\bar{1} \in C_{a_i} \subseteq \coker f$ and let $x_i\in X$ be determined by $f(x_i)=\pi^{a_i}y_i$. Then
\begin{align*}
    \tau^-_{\mathcal{S}(R)}f \cong
\begin{pmatrix}
-x_1 & \cdots & -x_r\\
\pi^{a_1} & & 0\\
& \ddots &\\
0 & & \pi^{a_r}
\end{pmatrix}
\colon \bigoplus_{i=1}^r C_{n-a_i} \longrightarrow X\oplus R^r
\end{align*}
where $x_i$ is identified with the morphism $C_{n-a_i}\to X, \bar{1} \mapsto x_i$. In particular $\tau^-_{\mathcal{S}(R)} f$ is contained in $\mathcal{S}(R')$ if and only if $f$ is an isomorphism.
\end{lem}
\begin{proof} Let $g$ denote the map, which we want to show coincides with $\smash{\tau^-_{\mathcal{S}(R)} f}$. The elements $y_i \in Y$ induce a map $p\colon R^r \to Y$ and we prove that
\begin{align*}
    0\longrightarrow \bigoplus_{i=1}^r C_{n-a_i}
\xlongrightarrow{g}
X\oplus R^r
\xrightarrow{(f\,\,p)}
Y\longrightarrow 0
\end{align*}
is exact. Clearly $g$ is a monomorphism, and $(f\,\,p)$ is an epimorphism since the $y_i$ generate inside $\coker f$ by definition. The equality $f(x_i) = \pi^{a_i} y_i$ shows that the composition is zero. Counting lengths
\begin{align*}
    \mathrm{len}_R\,Y +\sum_{i=1}^r \mathrm{len}_R\,C_{n-a_i} = (\mathrm{len}_R\,Y - \sum_{i=1}^r \mathrm{len}_R\,C_{a_i})+rn = \mathrm{len}_R\,X+\mathrm{len}_R\,R^r
\end{align*}
implies exactness. By \cite[Theorem~5.1]{RS3} it follows that $\smash{\tau_{\mathcal{S}(R)}\,g \cong \mathrm{Mimo}\,\tau_R\,(f\,\,p)}$.
Here $\tau_R$ must be applied to a morphism as a functor between the stable module categories. In our case $\tau_R$ acts as the identity on the stable level and $(f\,\,p) \cong f$ inside the stable category. Lastly $\mathrm{Mimo}$ acts as the identity on monomorphisms and so $\tau_{\mathcal{S}(R)}\,g  \cong f$. 

To deduce $g \cong \smash{\tau^-_{\mathcal{S}(R)}}f$ it is left to check that $g$ contains no summand isomorphic to $1_R \colon R \to R$ or $0 \to R$. Excluding the first is trivial as $C_{n-a_i} \not \cong R$, and for the second observe that $R$ would be a summand of $\coker g \cong Y$ which is not possible since $f \in \mathcal{S}(R')$.
\end{proof}
\noindent
\textbf{Computations.} The results in this section were used to reduce the computations C1--C3 to matrix equations over $R$, which were solved by GPT-6 Astra (xhigh) and verified by the AI in Lean 4 using mathlib. We illustrate the case $n =4$.

Let $R'=R/(\pi^3)$ and suppose that the indecomposable objects in $\mathcal{S}(R')$ are given by the depiction in \Cref{E: C3}, that is
\begin{align*}
    C_i\xrightarrow{\pi^{j-i}} C_j \quad (j>0\leq i \leq j\leq 3), \qquad C_2\xrightarrow{\left(\begin{smallmatrix}
        \pi\\1
    \end{smallmatrix}\right)}C_3\oplus C_1
\end{align*}
are all of them, up to isomorphism. We can now use \Cref{L: tau} to determine the necessary Auslander-Reiten translates for C1. 

\begin{equation*}
\begin{array}{c|c}
X & \tau^-_{\mathcal{S}(R)}X \\ \hline
0 \to C_3 & C_1\xrightarrow{\pi^3}C_4\\
C_1 \xrightarrow{\pi^2} C_3 & C_2\xrightarrow{\left(\begin{smallmatrix}
    \pi^2 \\1
\end{smallmatrix}\right)}C_4\oplus C_1\\
C_2 \xrightarrow{\left(\begin{smallmatrix}
    \pi\\
    1
\end{smallmatrix}\right)} C_3 \oplus C_1     & C_2\xrightarrow{\left(\begin{smallmatrix}
    \pi^2 \\ \pi
\end{smallmatrix}\right)}C_4\oplus C_2\\

0 \to C_2 & C_2\xrightarrow{\pi^2}C_4\\
C_2 \xrightarrow{\pi} C_3 & C_3\xrightarrow{\left(\begin{smallmatrix}
    \pi \\1
\end{smallmatrix}\right)}C_4\oplus C_2\\
C_1 \xrightarrow{\pi} C_2 & C_3\xrightarrow{\left(\begin{smallmatrix}
    \pi \\1
\end{smallmatrix}\right)}C_4\oplus C_1\\
0 \to C_1 & C_3\xrightarrow{\pi}C_4
\end{array}
\end{equation*}

Negative signs arising in the formula of \Cref{L: tau} have been removed by simple automorphisms to match the description from \Cref{E: C3}. Note that the AI took \Cref{L: tau} as the definition of the translates and it did not use Auslander-Reiten theory (in particular, also not for the Lean verification). Let $Y_2, \dots, Y_8$ denote the objects $\smash{\tau_{\mathcal{S}(R)}^- X}$ in the list and set $Y_1 = 0 \to C_4$ as well as  $Y_9 =1\colon  C_4 \to C_4$. This finishes the computation C1.

For the computation C2 we must realize $E=
    \operatorname{End}_{\mathcal{S}(R)/
    \langle\mathcal{S}(R')\rangle}
    \left(\oplus_i Y_i\right)$
as a bound path algebra. This is achieved by giving a morphism of $k$-algebras $f\colon kQ \to E$ for a quiver $Q$ and relations $I \subseteq kQ$ such that $f$ is surjective, $f$ maps $I$ to zero, and $\dim kQ/I = \dim E$. Clearly, then $kQ/I \cong E$. For our case, the quiver and relations are given by 
\begin{equation*}
\begin{tikzcd}[row sep=small, column sep=small]
  & & &  &  5 \arrow[dr] &  & & \\
1 \arrow[r] & 2 \arrow[r]
  & 3 \arrow[r]
  & 4 \arrow[ur] \arrow[dr] & \circlearrowleft & 7 \arrow[r]
  & 8 \arrow[r] & 9. \\
  & & &  & 6 \arrow[ur]  &  & &
\end{tikzcd}
\end{equation*}
The commutativity relation of the square is indicated by the circular arrow. Note that there is no dependency on $R$. The AI was first used to guess the correct map $f$ by investigating finite examples such as $R = \mathbb{Z}/16\mathbb{Z}$ or $R= \mathbb{F}_3[x]/(x^4)$. Then, the hardest task was to determine the spaces $\mathrm{Hom}_{\mathcal{S}(R)/\langle \mathcal{S}(R') \rangle} (Y_i, Y_j)$ explicitly and formally verify the necessary computations. Once this has been achieved, checking that $f$ is surjective, $I$ is mapped to zero, and $\dim kQ/I = \dim E$ is rather simple. 

We illustrate the computation of $\mathrm{Hom}_{\mathcal{S}(R)/\langle \mathcal{S}(R') \rangle} (Y_4, Y_7)$. The vertex $i$ in the quiver corresponds to the module $Y_i$, so we expect the space to be 1-dimensional. An arbitrary morphism is given by a commutative diagram
\begin{equation*}
\begin{tikzcd}
    C_2 \arrow[d, "r", swap] \arrow[r, "\left(\begin{smallmatrix}
    \pi^2 \\ \pi
\end{smallmatrix}\right)"] & C_4\oplus C_2 \arrow[d, "\left( \begin{smallmatrix}
     a \,\, b\\
     c \,\, d
\end{smallmatrix}\right)"]\\
    C_3 \arrow[r, "\left(\begin{smallmatrix}
    \pi \\ 1
\end{smallmatrix}\right)"] & C_4\oplus C_1
\end{tikzcd}
\end{equation*}
and the assignment $\left(\begin{smallmatrix}
    a & b \\ c & d
\end{smallmatrix} \right) \mapsto \bar{a}\in R/(\pi)$ provides the desired identification with $k$. Indeed, the kernel consists precisely of those $\left(\begin{smallmatrix}
    a & b \\ c & d
\end{smallmatrix} \right)$ that are annihilated by $\pi^3$ which corresponds to all morphisms factoring through $\mathcal{S}(R')$ by \Cref{L: C2}. Moreover, choosing $r = \pi$ and $\left(\begin{smallmatrix}
    a & b \\ c & d
\end{smallmatrix} \right) = \left(\begin{smallmatrix}
    1 & 0 \\ 0 & 1
\end{smallmatrix} \right)$ shows that the assignment is surjective. Similar computations have been carried out by the AI for the computation of Hom-spaces and a Lean verification was provided.  \Cref{L: C2} was simply used as a definition for the AI for the space of morphisms factoring through $\mathcal{S}(R')$, so no formalization of ideal quotients was necessary. 

For C3, we must list $s_4-s_3 = 10$ non-isomorphic objects in $\mathcal{S}(R)\setminus \mathcal{S}(R')$ and compute that each endomorphism ring is local in $\mathcal{S}(R)$. With $Y_i$ we already listed 9 of them. The remaining one is extracted from \Cref{E: C3} as 
\begin{align*}
    C_1 \oplus C_3 \xrightarrow{\left(\begin{smallmatrix}
        \pi & 1 \\
        0 & \pi
    \end{smallmatrix}\right)} C_2 \oplus C_4.
\end{align*}
A similar computation as before shows that the endomorphism ring equals $k$ inside the ideal quotient and endomorphisms factoring through $\mathcal{S}(R')$ are nilpotent. The objects will automatically be non-isomorphic as already the underlying $R$-modules are non-isomorphic. 

\begin{rem}\label{R: final rem} \rm The computations show that for $n\leq 5$, up to isomorphism, all of the indecomposable objects in $\mathcal{S}(R)$ are given by the figures of height ${\leq}n$ in the depiction from \Cref{E: C3}. Thus, their description is independent of $R$.

The following shows all of the quivers $Q$ and respective relations $I$ determined in the computations for $2\leq n \leq 6$ in increasing order (visualized by the AI). Recall that $B = (kQ/I)^\op$ is the desired finite dimensional algebra in \Cref{T: main indep}.
\begin{equation*}
\hspace*{-5mm}
\begin{gathered}
\phantom{X}\\[-1em]
\begin{tikzpicture}[quiver]
  \node[vertex] (a) at (0,0) {$\bullet$};
  \node[vertex] (b) at (1,0) {$\bullet$};
  \node[vertex] (c) at (2,0) {$\bullet$};
  \draw[arr] (a) -- (b);
  \draw[arr] (b) -- (c);
\end{tikzpicture}
\\[1em]
\begin{tikzpicture}[quiver]
  \foreach \i in {1,...,5}{%
    \node[vertex] (v\i) at (\i,0) {$\bullet$};
  }
  \draw[arr] (v1) -- (v2);
  \draw[arr] (v2) -- (v3);
  \draw[arr] (v3) -- (v4);
  \draw[arr] (v4) -- (v5);
\end{tikzpicture}
\\[1em]
\begin{tikzpicture}[quiver]
  \node[vertex] (v1) at (0,0) {$\bullet$};
  \node[vertex] (v2) at (1,0) {$\bullet$};
  \node[vertex] (v3) at (2,0) {$\bullet$};
  \node[vertex] (v4) at (3,0) {$\bullet$};
  \node[vertex] (v5) at (4,.7) {$\bullet$};
  \node[vertex] (v6) at (4,-.7) {$\bullet$};
  \node[vertex] (v7) at (5,0) {$\bullet$};
  \node[vertex] (v8) at (6,0) {$\bullet$};
  \node[vertex] (v9) at (7,0) {$\bullet$};
  \draw[arr] (v1) -- (v2);
  \draw[arr] (v2) -- (v3);
  \draw[arr] (v3) -- (v4);
  \draw[arr] (v4) -- (v5);
  \draw[arr] (v4) -- (v6);
  \draw[arr] (v5) -- (v7);
  \draw[arr] (v6) -- (v7);
  \draw[arr] (v7) -- (v8);
  \draw[arr] (v8) -- (v9);
  \commmark{4}{0}
\end{tikzpicture}
\\[1em]
\begin{tikzpicture}[quiver]
  \node[vertex] (v1) at (-4,0) {$\bullet$};
  \node[vertex] (v2) at (-3,0) {$\bullet$};
  \node[vertex] (v3) at (-2,0) {$\bullet$};
  \node[vertex] (v4) at (-1,0) {$\bullet$};
  \node[vertex] (v5) at (0,0) {$\bullet$};
  \node[vertex] (v8) at (1,0) {$\bullet$};
  \node[vertex] (v11) at (2,0) {$\bullet$};
  \node[vertex] (v14) at (3,0) {$\bullet$};
  \node[vertex] (v6) at (-1.56,.9) {$\bullet$};
  \node[vertex] (v9) at (-1.56,2.7) {$\bullet$};
  \node[vertex] (v12) at (0,3.6) {$\bullet$};
  \node[vertex] (v15) at (1.56,2.7) {$\bullet$};
  \node[vertex] (v17) at (1.56,.9) {$\bullet$};
  \node[vertex] (v7) at (-.78,1.35) {$\bullet$};
  \node[vertex] (v10) at (-.78,2.25) {$\bullet$};
  \node[vertex] (v13) at (0,2.7) {$\bullet$};
  \node[vertex] (v16) at (.78,2.25) {$\bullet$};
  \node[vertex] (v18) at (.78,1.35) {$\bullet$};
  \foreach \source/\target in {
    1/2,2/3,3/4,4/5,5/8,8/11,11/14,
    5/6,5/7,6/9,6/10,7/10,9/12,10/12,10/13,
    12/15,13/15,13/16,15/17,15/18,16/18,17/5,18/5}{%
    \draw[arr] (v\source) -- (v\target);
  }
  \commmark{-1.10}{1.05}
  \commmark{-1.10}{2.55}
  \commmark{.15}{3.08}
  \commmark{1.07}{2.42}
  \commmark{1.13}{1.17}
  \draw[zero]
  (.27,.276) -- (0,.12) -- (-.27,.276);
\draw[zero]
  (.27,.416) -- (0,.26) -- (-.11,.450);
\draw[zero]
  (.11,.610) -- (0,.42) -- (-.27,.576);
\draw[zero]
  (.11,.770) -- (0,.58) -- (-.11,.770);
\draw[zero]
  (-.30,-.14) -- (.30,-.14);
\end{tikzpicture}
\\[2em]
\begin{tikzpicture}[quiver,x=4mm,y=4.8mm]
\draw[longzero,quintic]
  (19.6,-3.15) -- (22.05,-2.1) -- (24.5,-1.05)
  -- (26.95,0) -- (29.4,1.05) -- (31.85,2.1);
\fill[quintic] (19.6,-3.15) circle[radius=2.4pt];
\fill[quintic] (31.85,2.1) circle[radius=2.4pt];
\draw[longzero,quintic]
  (-2.45,2.1) -- (0,1.05) -- (2.45,0)
  -- (4.9,-1.05) -- (7.35,-2.1) -- (9.8,-3.15);
\fill[quintic] (-2.45,2.1) circle[radius=2.4pt];
\fill[quintic] (9.8,-3.15) circle[radius=2.4pt];
%
\draw[longzero,cubic]
  (4.9,-3.15) -- (7.35,-2.1)
  -- (9.8,-1.05) -- (12.25,-1.05);
\fill[cubic] (4.9,-3.15) circle[radius=2.4pt];
\fill[cubic] (12.25,-1.05) circle[radius=2.4pt];
\draw[longzero,cubic]
  (-2.45,-1.05) -- (0,-1.05) -- (2.45,0) -- (4.9,1.05);
\fill[cubic] (-2.45,-1.05) circle[radius=2.4pt];
\fill[cubic] (4.9,1.05) circle[radius=2.4pt];
\draw[longzero,cubic]
  (17.15,-1.05) -- (19.6,-1.05)
  -- (22.05,-2.1) -- (24.5,-3.15);
\fill[cubic] (17.15,-1.05) circle[radius=2.4pt];
\fill[cubic] (24.5,-3.15) circle[radius=2.4pt];
\draw[longzero,cubic]
  (24.5,1.05) -- (26.95,0)
  -- (29.4,-1.05) -- (31.85,-1.05);
\fill[cubic] (24.5,1.05) circle[radius=2.4pt];
\fill[cubic] (31.85,-1.05) circle[radius=2.4pt];
%
\draw[guide,line width=.45pt]
  (-.2,1.4) -- (-.3,1.4) -- (-.3,-1.4) -- (-.2,-1.4);
\draw[guide,line width=.45pt]
  (29.6,1.4) -- (29.7,1.4) -- (29.7,-1.4) -- (29.6,-1.4);
%
\node[vertex,text=ghost] (g47) at (-2.45,-2.1) {$\bullet$};
\node[vertex,text=ghost] (g46) at (-2.45,-1.05) {$\bullet$};
\node[vertex,text=ghost] (g45) at (-2.45,0) {$\bullet$};
\node[vertex,text=ghost] (g44) at (-2.45,2.1) {$\bullet$};
\node[vertex] (v2) at (0,-1.05) {$\bullet$};
\node[vertex] (v1) at (0,1.05) {$\bullet$};
\node[vertex] (v6) at (2.45,-2.1) {$\bullet$};
\node[vertex] (v5) at (2.45,-1.05) {$\bullet$};
\node[vertex] (v4) at (2.45,0) {$\bullet$};
\node[vertex] (v3) at (2.45,2.1) {$\bullet$};
\node[vertex] (v11) at (4.9,-3.15) {$\bullet$};
\node[vertex] (v10) at (4.9,-1.05) {$\bullet$};
\node[vertex] (v9) at (4.9,1.05) {$\bullet$};
\node[vertex] (v8) at (4.9,3.15) {$\bullet$};
\node[vertex] (v7) at (4.9,5.25) {$\bullet$};
\node[vertex] (v16) at (7.35,-2.1) {$\bullet$};
\node[vertex] (v15) at (7.35,-1.05) {$\bullet$};
\node[vertex] (v14) at (7.35,0) {$\bullet$};
\node[vertex] (v13) at (7.35,2.1) {$\bullet$};
\node[vertex] (v12) at (7.35,4.2) {$\bullet$};
\node[vertex] (v20) at (9.8,-3.15) {$\bullet$};
\node[vertex] (v19) at (9.8,-1.05) {$\bullet$};
\node[vertex] (v18) at (9.8,1.05) {$\bullet$};
\node[vertex] (v17) at (9.8,3.15) {$\bullet$};
\node[vertex] (v23) at (12.25,-2.1) {$\bullet$};
\node[vertex] (v24) at (12.25,-1.05) {$\bullet$};
\node[vertex] (v22) at (12.25,0) {$\bullet$};
\node[vertex] (v21) at (12.25,2.1) {$\bullet$};
\node[vertex] (v26) at (14.7,-1.05) {$\bullet$};
\node[vertex] (v25) at (14.7,1.05) {$\bullet$};
\node[vertex] (v30) at (17.15,-2.1) {$\bullet$};
\node[vertex] (v29) at (17.15,-1.05) {$\bullet$};
\node[vertex] (v28) at (17.15,0) {$\bullet$};
\node[vertex] (v27) at (17.15,2.1) {$\bullet$};
\node[vertex] (v34) at (19.6,-3.15) {$\bullet$};
\node[vertex] (v33) at (19.6,-1.05) {$\bullet$};
\node[vertex] (v32) at (19.6,1.05) {$\bullet$};
\node[vertex] (v31) at (19.6,3.15) {$\bullet$};
\node[vertex] (v39) at (22.05,-2.1) {$\bullet$};
\node[vertex] (v38) at (22.05,-1.05) {$\bullet$};
\node[vertex] (v37) at (22.05,0) {$\bullet$};
\node[vertex] (v36) at (22.05,2.1) {$\bullet$};
\node[vertex] (v35) at (22.05,4.2) {$\bullet$};
\node[vertex] (v43) at (24.5,-3.15) {$\bullet$};
\node[vertex] (v42) at (24.5,-1.05) {$\bullet$};
\node[vertex] (v41) at (24.5,1.05) {$\bullet$};
\node[vertex] (v40) at (24.5,3.15) {$\bullet$};
\node[vertex] (v47) at (26.95,-2.1) {$\bullet$};
\node[vertex] (v46) at (26.95,-1.05) {$\bullet$};
\node[vertex] (v45) at (26.95,0) {$\bullet$};
\node[vertex] (v44) at (26.95,2.1) {$\bullet$};
\node[vertex] (v2r) at (29.4,-1.05) {$\bullet$};
\node[vertex] (v1r) at (29.4,1.05) {$\bullet$};
\node[vertex,text=ghost] (g5) at (31.85,-1.05) {$\bullet$};
\node[vertex,text=ghost] (g3) at (31.85,2.1) {$\bullet$};
%
\draw[arr] (v1) -- (v3);
\draw[arr] (v1) -- (v4);
\draw[arr] (v2) -- (v4);
\draw[arr] (v2) -- (v5);
\draw[arr] (v2) -- (v6);
\draw[arr] (v3) -- (v8);
\draw[arr] (v3) -- (v9);
\draw[arr] (v4) -- (v9);
\draw[arr] (v4) -- (v10);
\draw[arr] (v5) -- (v10);
\draw[arr] (v6) -- (v10);
\draw[arr] (v6) -- (v11);
\draw[arr] (v7) -- (v12);
\draw[arr] (v8) -- (v13);
\draw[arr] (v9) -- (v13);
\draw[arr] (v9) -- (v14);
\draw[arr] (v10) -- (v14);
\draw[arr] (v10) -- (v15);
\draw[arr] (v10) -- (v16);
\draw[arr] (v11) -- (v16);
\draw[arr] (v12) -- (v17);
\draw[arr] (v13) -- (v18);
\draw[arr] (v14) -- (v18);
\draw[arr] (v14) -- (v19);
\draw[arr] (v15) -- (v19);
\draw[arr] (v16) -- (v19);
\draw[arr] (v16) -- (v20);
\draw[arr] (v17) -- (v21);
\draw[arr] (v18) -- (v22);
\draw[arr] (v19) -- (v22);
\draw[arr] (v19) -- (v23);
\draw[arr] (v19) -- (v24);
\draw[arr] (v20) -- (v23);
\draw[arr] (v21) -- (v25);
\draw[arr] (v22) -- (v25);
\draw[arr] (v22) -- (v26);
\draw[arr] (v23) -- (v26);
\draw[arr] (v24) -- (v26);
\draw[arr] (v25) -- (v27);
\draw[arr] (v25) -- (v28);
\draw[arr] (v26) -- (v28);
\draw[arr] (v26) -- (v29);
\draw[arr] (v26) -- (v30);
\draw[arr] (v27) -- (v31);
\draw[arr] (v28) -- (v32);
\draw[arr] (v28) -- (v33);
\draw[arr] (v29) -- (v33);
\draw[arr] (v30) -- (v33);
\draw[arr] (v30) -- (v34);
\draw[arr] (v31) -- (v35);
\draw[arr] (v32) -- (v36);
\draw[arr] (v32) -- (v37);
\draw[arr] (v33) -- (v37);
\draw[arr] (v33) -- (v38);
\draw[arr] (v33) -- (v39);
\draw[arr] (v34) -- (v39);
\draw[arr] (v36) -- (v40);
\draw[arr] (v36) -- (v41);
\draw[arr] (v37) -- (v41);
\draw[arr] (v37) -- (v42);
\draw[arr] (v38) -- (v42);
\draw[arr] (v39) -- (v42);
\draw[arr] (v39) -- (v43);
\draw[arr] (v40) -- (v44);
\draw[arr] (v41) -- (v44);
\draw[arr] (v41) -- (v45);
\draw[arr] (v42) -- (v45);
\draw[arr] (v42) -- (v46);
\draw[arr] (v42) -- (v47);
\draw[arr] (v43) -- (v47);
\draw[arr] (v44) -- (v1r);
\draw[arr] (v45) -- (v1r);
\draw[arr] (v45) -- (v2r);
\draw[arr] (v46) -- (v2r);
\draw[arr] (v47) -- (v2r);
%
\draw[arr,draw=ghost] (v2r) -- (g5);
\draw[arr,draw=ghost] (g44) -- (v1);
\draw[arr,draw=ghost] (g45) -- (v1);
\draw[arr,draw=ghost] (g45) -- (v2);
\draw[arr,draw=ghost] (g46) -- (v2);
\draw[arr,draw=ghost] (g47) -- (v2);
\draw[arr,draw=ghost] (v1r) -- (g3);
%
\draw[zero]
  (6.8,-2.536) -- (7.35,-2.3) -- (7.9,-2.536);
\draw[zero]
  (21.5,-2.536) -- (22.05,-2.3) -- (22.6,-2.536);
\draw[zero]
  (-.6,-.87) -- (.6,-.87);
\draw[zero]
  (9.2,-.87) -- (10.4,-.87);
\draw[zero]
  (9.25,-1.25) -- (9.8,-1.25) -- (10.35,-1.486);
\draw[zero]
  (19,-.87) -- (20.2,-.87);
\draw[zero]
  (-.55,1.486) -- (0,1.25) -- (.55,1.486);
\draw[zero]
  (4.35,-1.036) -- (5.45,-.564);
\draw[zero]
  (11.7,.436) -- (12.25,.2) -- (12.8,.436);
\draw[zero]
  (14.15,.614) -- (14.7,.85) -- (15.25,.614);
\draw[zero]
  (14.15,-.614) -- (14.7,-.85) -- (15.25,-.614);
\draw[zero]
  (19.05,-1.486) -- (19.6,-1.25) -- (20.15,-1.25);
\draw[zero]
  (23.95,-1.064) -- (25.05,-1.536);
\draw[zero]
  (14.15,1.486) -- (14.7,1.25) -- (15.25,1.486);
\draw[zero]
  (4.3,-1.23) -- (5.5,-1.23);
\draw[zero]
  (16.6,.436) -- (17.15,.2) -- (17.7,.436);
\draw[zero]
  (23.9,-.87) -- (25.1,-.87);
\draw[zero]
  (-.55,-1.486) -- (0,-1.25) -- (.55,-1.486);
\draw[zero]
  (14.1,-.99) -- (15.3,-.99);
\draw[zero]
  (14.15,-1.486) -- (14.7,-1.25) -- (15.25,-1.486);
  %
\commmark{4.9}{-2.1}
\commmark{26.95}{1.05}
\commmark{19.6}{-2.1}
\commmark{24.5}{0}
\commmark{24.5}{2.1}
\commmark{9.8}{-2.1}
\commmark{22.05}{1.05}
\commmark{24.5}{-2.1}
\commmark{2.45}{1.05}
\commmark{4.9}{0}
\commmark{4.9}{2.1}
\commmark{7.35}{1.05}
\commmark{9.8}{0}
\commmark{0}{0}
\commmark{19.6}{0}
%
\summark{2.45}{-1.05}
\summark{7.35}{-1.05}
\summark{26.95}{-1.05}
\summark{12.25}{-1.05}
\summark{22.05}{-1.05}
\summark{17.15}{-1.05}
\end{tikzpicture}\\[1em]
\end{gathered}
\end{equation*}

Dotted lines indicate that the corresponding compositions of consecutive
arrows vanish, $\circlearrowleft$ indicates that the surrounding square commutes, and $\vcenter{\hbox{\begin{tikzpicture}
  \summark{0}{0}
\end{tikzpicture}}}$ indicates that the composition of the two arrows through the marked vertex equals the sum of the two outer paths with the same endpoints. For $n=6$ the bracketed ends are identified. The green and blue highlighted paths represent zero relations.
\end{rem}

\end{document}